\documentclass[11pt]{article}
\usepackage{amssymb}
\usepackage{subfigure}
\usepackage{amsbsy}
\usepackage{amsmath}
\usepackage{graphicx}
\usepackage[varg]{pxfonts}
\usepackage{epstopdf}
\usepackage{color}

\renewcommand {\thefootnote}{\fnsymbol{footnote}}

\numberwithin{equation}{section} \allowdisplaybreaks

\begin{document}
\title{\Large \bf A Gauss curvature flow related to the Orlicz-Aleksandrov problem
\thanks{Research is supported in part by the Natural Science
Foundation of China (No.11871275; No.11371194).  }}
\author{ \small \bf Bin Chen$^{1}$\thanks{E-mail: chenb121223@163.com\ bchen@njust.edu.cn.}, and \bf Peibiao Zhao$^{1}$\thanks{Corresponding
author E-mail: pbzhao@njust.edu.cn.}\\
\\
\small      (1. Department of Mathematics, Nanjing
University of Science and Technology, Nanjing, China)
}
\date{}
\maketitle
\renewcommand{\thefootnote}{\fnsymbol{footnote}}

\vskip 20pt

\begin{center}
\begin{minipage}{12.8cm}
\small
 {\bf Abstract:}  In this paper we first obtain the existence of smooth solutions to Orlicz-Aleksandrov problem via a Gauss-like curvature flow.



 \vskip 0.1cm
\end{minipage}
\end{center}

\vskip 20pt
\section{\bf Introduction}
The curvature measure of convex bodies is one of the basic
principles in convex geometry analysis. In particular, it plays key
role in the Brunn-Minkowski theory of convex bodies. The most
studied of the curvature measures is the Aleksandrov's integral
curvature (also called integral Gauss curvature) defined and studied
by Aleksandrov \cite{Al} using a topological argument.
Moreover, the famous Aleksandrov problem
with respect to integral curvatures is an important cornerstone of
the Brunn-Minkowski theory.

The integral curvature, $\mathcal{J}(K,\cdot)$ on the unit sphere
$S^{n-1}$, of $K\in\mathcal{K}_o^n$ (the set of convex bodies
containing the origin in their interiors) is defined by
$$\mathcal{J}(K,\omega)=\mathcal{H}^{n-1}(\Re_K(\omega))$$
for each Borel set $\omega\subset S^{n-1}$. Here $\mathcal{H}^{n-1}$ is the $(n-1)$-dimensional Hausdorff measure, and $\Re_K$ is the radial Gauss image.

The classical Aleksandrov problem, roughly speaking, asks for
necessary and sufficient conditions for a given Borel measure $\mu$
on the unit sphere so that the measure is the integral curvature of
a convex body in $\mathbb{R}^n$. Namely, this problem is to find a
convex body $K\subset\mathbb{R}^n$ such that
\begin{align}\label{1.0}
d\mathcal{J}(K,\cdot)=d\mu\ \ on \ \ S^{n-1}.
\end{align}

The PDE associated with classical Aleksandrov problem for integral
curvature asks (see \cite{O} or \cite{O1}): If the given measure
$\mu$ has a density $g$, then the (\ref{1.0}) is equivalent to
solving the following Monge-Amp\`{e}re type equation
$$\frac{h}{(|\nabla h|^2+h^2)^{\frac{n}{2}}}\det(\nabla^2h+hI)=g,$$
where $h$ is the support function of convex body $K$, $\nabla h$ is the gradient of $h$, while $\nabla^2h$ is the
Hessian matrix of $h$ with respect to an orthonormal frame on $S^{n-1}$, and $I$ is the identity matrix of order $n-1$.

In the Huang et al's work \cite{HL}, for all real $p$, there is a
geometrically natural $L_p$ extension of integral curvature (also
called $L_p$-integral curvature). To state the $L_p$-Aleksandrov
problem with respect to $L_p$-integral curvature, in full
generality, it is necessary to introduce the entropy functional
$$\mathcal{E}(K)=-\int_{S^{n-1}}\log h_K(v)dv$$
for $K\in\mathcal{K}_o^n$. Here the integration is with respect to spherical Lebesgue measure.

For $p\neq0$, the $L_p$-integral curvature, $\mathcal{J}_p(K,\cdot)$, of $K\in \mathcal{K}_o^n$, as a Borel measure on $S^{n-1}$ is defined by the variational formula
\begin{align}\label{1.1}
\frac{d}{dt}\mathcal{E}(K\hat{+}_pt\cdot Q)\bigg|_{t=0}=\frac{1}{p}\int_{S^{n-1}}\rho_Q^{-p}(u)d\mathcal{J}_p(K,u)
\end{align}
which holds for each $Q\in \mathcal{K}_o^n$, where $\rho_Q$ is the radial function and $K\hat{+}_pt\cdot Q$ is the harmonic $L_p$-combination (see \cite{HL}). It turns out that for each $K\in \mathcal{K}_o^n$,
$$d\mathcal{J}_p(K,\cdot)=\rho_K^pd\mathcal{J}(K,\cdot).$$
When $p=0$ in the harmonic $L_p$-combination, the integral curvature, $\mathcal{J}(K,\cdot)$, of $K\in \mathcal{K}_o^n$ is a Borel measure on $S^{n-1}$ that can be defined by
\begin{align}\label{1.2}
\frac{d}{dt}\mathcal{E}(K\hat{+}_0t\cdot Q)\bigg|_{t=0}=-\int_{S^{n-1}}\log\rho_Q(u)d\mathcal{J}(K,u),
\end{align}
which holds for $Q\in \mathcal{K}_o^n$. It should be emphasized here
that the variational formula (\ref{1.2}) is not Aleksandrov's
definition of classical integral curvatures.

\noindent{\bf The $L_p$-Aleksandrov problem}~~{\it For a fixed
$p\in\mathbb{R}$, and a given Borel measure $\mu$ on $S^{n-1}$, what are the necessary and sufficient conditions
so that
\begin{align}\label{1.2.0}
d\mathcal{J}_p(K,\cdot)=d\mu
\end{align}
of a convex body $K\in\mathcal{K}_o^n$?}

Moreover, the PDE associated with the $L_p$-Aleksandrov problem asks
(see \cite{HL}): If the given measure $\mu$ has a density $g$, then
the (\ref{1.2.0}) is equivalent to solving the following
Monge-Amp\`{e}re type equation
\begin{align}\label{1.2.01}
\frac{h^{1-p}}{(|\nabla h|^2+h^2)^{\frac{n}{2}}}\det(\nabla^2h+hI)=g.
\end{align}

{For the $L_p$ Aleksandrov problem (\ref{1.2.0}), when $p>0$, the
existence of solutions (measure solutions) has been completely
solved (see \cite{HL}, Theorem 7.1); When $p=0$, the $L_p$
Aleksandrov problem (\ref{1.2.0}) is exactly classical Aleksandrov
problem;
  When $p<0$ and the measure is even, the
sufficient condition for the existence of the solution to the $L_p$
Aleksandrov problem (\ref{1.2.0}) is given (see \cite{HL}, Theorem
7.3).  However, there are still many importmant problems, such as
non-even solution, smoothness of the solution, etc., has not been
solved.}

Recently, the concept of the Orlicz-integral curvature
$\mathcal{J}_{\phi}(K,\cdot)$ of convex body $K\in\mathcal{K}_o^n$
was defined by the following variational formula (see \cite{FH1} in
detail): Let $\Omega\subset S^{n-1}$ be a closed set not contained
in any closed hemisphere of $S^{n-1}$. Denote  $C(\Omega)$ by the
set of continuous functions on $\Omega$. For $f\in C(\Omega)$ and
$\rho_{K}\in C^+(\Omega)$. If $\phi$ is of continuously
differentiable and strictly monotonic function on $(0,\infty)$, then

$$\frac{d}{dt}\mathcal{E}(\langle \rho_t\rangle)\bigg|_{t=0}
=\int_{\Omega}f(u)d\mathcal{J}_\phi(K,u),$$ where
$\rho_t(u)=\phi^{-1}(\phi(\rho_K(u))+tf(u))$ and $\langle
\rho_t\rangle$ denotes the convex hull. It turns out that for
$K\in\mathcal{K}_o^n$
\begin{align}\label{1.2.1}
d\mathcal{J}_\phi(K,\cdot)=\phi(\rho_K)d\mathcal{J}(K,\cdot).
\end{align}
When $\phi(t)=t^{p}$ in (\ref{1.2.1}), it is just the $L_p$ integral curvature.

The following Orlicz-Aleksandrov problem was proposed in \cite{FH1}:

\noindent{\bf The Orlicz-Aleksandrov problem}~~ {\it For a suitable
continuous function $\phi: (0,\infty)\rightarrow(0,\infty)$, and  a
non-zero finite Borel measure $\mu$ on $S^{n-1}$, do there exists a
constant $\lambda>0$ and a convex body $K\in\mathcal{K}_o^n$ such
that
\begin{align}\label{1.3}
\lambda d\mathcal{J}_{\phi}(K,\cdot)=d\mu\ \ on \ \ S^{n-1}?
\end{align}}
For the Orlicz-Aleksandrov problem, when the given measure is even, was first solved in two situations via the variational method, see, e.g., \cite{FH1}.

We note that when the given measure $\mu$ has  a density $g$, then the
Orlicz-Aleksandrov problem (\ref{1.3}) is equivalent to solving the
following  Monge-Amp\`{e}re type equation (see Section 3 in detail):
\begin{align}\label{1.4}
\frac{\lambda h\phi(1/h)}{(|\nabla h|^2+h^2)^{\frac{n}{2}}}\det(\nabla^2h+hI)=g.
\end{align}

We know that the Orlicz Aleksandrov problem is a generalization of the classical Aleksandrov problem \cite{Al,O,O1}. When $\phi(t)=t^p$ with $t\in\mathbb{R}$, Eq. (\ref{1.4}) corresponds to the $L_p$ Aleksandrov problem \cite{HL,Z}.

{Moreover, the smoothness of solutions to Aleksandrov type problems
(or Minkowski type problems) and the non-even smooth solution of the
Orlicz Aleksandrov problems (as well as related Monge-Amp\`{e}re
equation) are open.}

Recently, the argument of   the smoothness of the even-solutions of
Minkowski type problems via the geometric flow method has been made
great progress(see  \cite{B1,CH,CW,LL,LS} for details).

Motivated by the above statements,  we first in this paper study the
existence of smooth non-even solution to the Eq. (\ref{1.4}) with
$\lambda=1$. The following theorem shows the existence of the smooth
solution to the Orlicz-Aleksandrov problem.

\noindent{\bf Theorem 1.1}~~{\it Suppose $\phi: (0,+\infty)\rightarrow(0,+\infty)$ is a continuous function. For any given positive smooth function $g$ on $S^{n-1}$ satisfying
\begin{align}\label{1.5}
\limsup_{s\rightarrow+\infty}\{\phi(s)\}<g<\liminf_{s\rightarrow0^+}\{\phi(s)\},
\end{align}
then the equation (\ref{1.4}) has a smooth solution $h$ with $\lambda=1$.}

In order to obtain the Theorem 1.1, our main idea is reflected in
the following two folds: I). Find a suitable anisotropic Gauss-like
curvature flow; II). Find a monotone functional of the solution to
the flow, which is the key to prove the existence of a solution to
Eq.(\ref{1.4}).

Let $M_0$ be a closed, smooth and strictly convex hypersurface in $\mathbb{R}^n$ enclosing the origin and given by a smooth
embedding $X_0: S^{n-1}\rightarrow\mathbb{R}^n$. In this paper we consider a family of closed hypersurfaces $M_t$ given by smooth maps $X: S^{n-1}\times [0,T)\rightarrow\mathbb{R}^n$ satisfying the initial value problem:
\begin{align}\label{1.6}
\begin{cases}
\frac{\partial X(x,t)}{\partial t}=-g(\nu)\frac{r^n}{\phi(r)}\mathcal{K}\nu+X(x,t),\\
X(x,0)=X_0(x),
\end{cases}
\end{align}
where $g$ is a given positive smooth function on $S^{n-1}$, $r=|X|$ is the distance from $X$ to the origin, $\phi: (0,+\infty)\rightarrow(0,+\infty)$ is a positive smooth function, $\mathcal{K}$ is the Gauss curvature of the hypersurface $M_t$ parametrized by $X(x,t)$,  $\nu$ is the unit outer normal vector at $X(x,t)$, and $T$ is the maximal time for which the solution exists.

The Gauss curvature flow was introduced by Firey \cite{F} to model the shape change of worn stones.
Since then,
various Gauss curvature flows have been extensively studied, see,
e.g. \cite{A,A0,AG,B,B1,CH,CW,G0,G2,I1,I,LL,LS,Wx}. In addition, the method of geometric flow
to solve some famous geometric inequalities
has also attracted the attention of many scholars, see, e.g. \cite{AC,G1,H1,Hl,HLW,Ly,LWX,WX}.

We will show that the flow (\ref{1.6}) has a long-time solution, and
derive that the support function of limiting hypersurface of this
flow provides a smooth solution to Eq. (\ref{1.4}). The following
functional related to the flow (\ref{1.6}) plays an important role
in our argument,
$$\mathcal{F}(M_t)=\int_{S^{n-1}}\log h(x,t)dx-\int_{S^{n-1}}\frac{\varphi(r(\xi,t))}{g(x)}d\xi,\ \ \ \ (\cdot,t)\in S^{n-1}\times[0,T),$$
where $h$, $r$ are the support function and radial function of $M_t$, and
$$\varphi(t)=\int_0^t\frac{\phi(s)}{s}ds.$$

Now, we obtain the long-time existence and convergence of the flow (\ref{1.6}).

\noindent{\bf Theorem 1.2}~~{\it Let $M_0$ be a closed, smooth, and uniformly convex hypersurface in $\mathbb{R}^n$. Assume functions $g$ and $\phi$ satisfy the assumptions of Theorem 1.1, then the flow (\ref{1.6}) has a smooth solution $M_t$, which exists for any time $t\in[0,\infty)$.
{Moreover, when $t\rightarrow\infty$, the support function $h_t$ of $M_t$ converges in $C^\infty$ to a smooth solution $h_\infty$ to (\ref{1.4}) with $\lambda=1$, which is the support function of a smooth, closed, and uniformly convex hypersurface $M_\infty$.}}

The organization is as follows. The corresponding background materials and some results are introduced in Section 2. In Section 3, we establish the Gauss-like curvature flow and related functional. In Section 4,
we obtain the long-time existence of the flow (\ref{1.6}). In section 5, we prove the Theorem 1.2, and provide a special uniqueness of Orlicz-Aleksandrov problem.

\section{\bf Preliminaries}

{\bf 2.1~~Convex body and Orlicz integral curvature}

The basic facts in this subsection can be found in Gardner and Schneider's book \cite{G,S}, which are the standard references regarding convex bodies, and references \cite{H,HL}.
Let $\mathbb{R}^n$ denote the $n$-dimensional Euclidean space. The unit sphere in $\mathbb{R}^n$ is denoted by $S^{n-1}$. A convex body in $\mathbb{R}^n$ is a compact convex set with nonempty interior. Denote by $\mathcal{K}_o^n$  the class of convex bodies in $\mathbb{R}^n$ that contain the origin in their interiors. Let $K\in\mathcal{K}_o^n$, the radial function $\rho_K: \mathbb{R}^n\backslash\{0\}\rightarrow\mathbb{R}$ is defined by
$$\rho_K(x)=\max\{\lambda: \lambda x\in K\}, \ \ \ \ x\in\mathbb{R}^n\backslash\{0\}.$$
For $u\in S^{n-1}$, there is $\rho_K(u)u\in\partial K$.

The support function, $h_K: S^{n-1}\rightarrow\mathbb{R}$, of a convex body $K$ in $\mathbb{R}^n$ is defined by
$$h_K(u)=\max\{u\cdot x: x\in K\}$$
for $u\in S^{n-1}$, where $u\cdot x$ is the standard inner product of $u$ and $x$ in $\mathbb{R}^n$.

The radial function and the support function are related,
$$h_K(v)=\sup_{u\in S^{n-1}}\{\rho_K(u)u\cdot v\},$$
$$\frac{1}{\rho_K(u)}=\sup_{v\in S^{n-1}}\frac{u\cdot v}{h_K(v)}.$$

For a convex body $K\in\mathcal{K}_o^n$, the polar body $K^*$ of $K$ is
$$K^*=\{x\in\mathbb{R}^n: x\cdot y\leq 1, \ \ for\ \ all\ \ y\in K\}.$$
The support function and radial function of the convex body and its polar are related in the following way,
\begin{align}\label{2.1.1}
h_K(x)=\frac{1}{\rho_{K^*}(x)},\ \ \ \ \rho_K(x)=\frac{1}{h_{K^*}(x)}.
\end{align}

The integral curvature, $\mathcal{J}(K,\cdot)$, of $K\in\mathcal{K}_o^n$ is a Borel measure on $S^{n-1}$ defined by
\begin{align*}
\mathcal{J}(K,\omega)=\mathcal{H}(\Re_K(\omega))
\end{align*}
for each Borel set $\omega\subset S^{n-1}$, where radial Gauss image $\Re_K(\omega)$ of $\omega$ given by
$$\Re_K(\omega)=\{u\in S^{n-1}:\rho_K(v)v\in H_K(u)\ \ for\ \ some\ \ v\in\omega\},$$
where $H_K$ is the supporting hyperplane of $K$ with the outer unit normal $u$.

If the convex body $K$ is $C^2$ smooth with positive Gauss curvature, then the integral curvature has a continuous
density (see, e.g., \cite{HL,O}),
\begin{align}\label{2.1.2}
\frac{h}{(|\nabla h|^2+h^2)^{\frac{n}{2}}}\det(\nabla^2h+hI),
\end{align}
where $h=\rho^{-1}_K$, while $\nabla h$ and $\nabla^2h$ are the gradient and the Hessian matrix of $h$, and $I$ is the identity matrix, with respect to an orthonormal frame on $S^{n-1}$.

For $K\in\mathcal{K}_o^n$ and $\phi: (0,\infty)\rightarrow(0,\infty)$ is a continuous function. The Orlicz-integral curvature was defined by (see \cite{FH1})
$$\mathcal{J}_\phi(K,\omega)=\int_{\Re_K(\omega)}\phi(\rho_K(\alpha_K^*(u)))du$$
for each Borel set $\omega\subset S^{n-1}$. Here $\alpha_K^*(u)$ is the reverse radial Gauss map of $K$.
Moreover, the Orlicz-integral curvature is absolutely continuous with respect to the classical integral curvature $\mathcal{J}(K,\cdot)$, namely
\begin{align}\label{2.1.3}
d\mathcal{J}_\phi(K,\cdot)=\phi(\rho_K)d\mathcal{J}(K,\cdot).
\end{align}
Obviously, when $\phi(t)=t^p$ ($p\in\mathbb{R}$), the Orlicz-integral curvature is just the $L_p$-integral curvature introduced by Huang et al  \cite{HL}
$$\mathcal{J}_p(K,\omega)=\int_{\Re_K(\omega)}\rho_K^p(\alpha_K^*(u))du.$$

\noindent{\bf 2.2~~Convex hypersurfaces}

Let $M$ be a closed, smooth, uniformly convex hypersurfaces in $\mathbb{R}^n$. Assume that $M$ is parametrized by the inverse Gauss map
\begin{align*}
X=\upsilon_M^{-1}: S^{n-1}\rightarrow M.
\end{align*}
The support function $h: S^{n-1}\rightarrow\mathbb{R}$ of $M$ is defined by
\begin{align*}
h(x)=\max\{\langle x,y\rangle,\ \ y\in M\},\ \ \ \ x\in S^{n-1}.
\end{align*}
The supremum is attained at a point $y$ such that $x$ is the outer normal of $M$ at $X$. It is easy to check that
\begin{align}\label{2.2.1}
X(x)=h(x)x+\nabla h(x),
\end{align}
where $\nabla$ is the covariant derivative with respect to the standard metric $\delta_{ij}$ of the sphere $S^{n-1}$.
Denote the radial function of $M_t$ by $\rho(u,t)$. From (\ref{2.2.1}), $u$ and $x$ are related by
\begin{align}\label{2.2.11}
\rho(u)u=h(x)x+\nabla h(x),
\end{align}
then there is
\begin{align*}
x=\frac{\rho(u)u-\nabla\rho}{\sqrt{\rho^2+|\nabla\rho|^2}}.
\end{align*}
From the definitions of radial function and $r$, then
\begin{align}\label{2.2.2}
r=|X|=\bigg(|\nabla h|^2+h^2\bigg)^{\frac{1}{2}}
\end{align}
and
\begin{align}\label{2.2.3}
h=\frac{r^2}{\sqrt{r^2+|\nabla r|^2}}.
\end{align}

The second fundamental form $A_{ij}$ of $M$ can be computed in terms of the support function (see e.g., \cite{A1,U})
\begin{align}\label{2.2.4}
A_{ij}=\nabla_{ij}h+he_{ij},
\end{align}
where $\nabla_{ij}=\nabla_i\nabla_j$ denotes the second order covariant derivative with respect to $e_{ij}$. The induced metric matix $g_{ij}$ of $M$ can be derived by Weingarten's formula,
\begin{align}\label{2.2.5}
e_{ij}=\langle\nabla_ix,\nabla_jx\rangle=A_{ik}A_{lj}g^{kl}.
\end{align}
It follows from (\ref{2.2.4}) and (\ref{2.2.5}) that the principal radii of curvature of $M$, under a smooth local orthonormal frame on $S^{n-1}$, are the eigenvalues of the matrix
\begin{align*}
b_{ij}=\nabla_{ij}h+h\delta_{ij}.
\end{align*}
We will use $b^{ij}$ to denote the inverse matrix of $b_{ij}$. In particular, the Gauss curvature is given by
\begin{align}\label{2.2.6}
\mathcal{K}(x)=(\det(\nabla_{ij}h+h\delta_{ij}))^{-1}=S_n^{-1}(\nabla_{ij}h+h\delta_{ij}),
\end{align}
where
\begin{align*}
S_k=\sum_{i_1<\cdot\cdot\cdot<i_k}\lambda_{i_1}\cdot\cdot\cdot\lambda_{i_k}
\end{align*}
denotes the $k$-th elementary symmetric polynomial.

Let $r$, $\alpha$ and $\alpha^*$ be the radial function, radial Gauss mapping and reverse radial Gauss mapping of $M$.
It is well-known that the determinants of the Jacobian of radial Gauss mapping $\alpha$ and reverse radial Gauss mapping of $M$ are given by, see e.g. \cite{H},
\begin{align}\label{2.2.7}
|Jac \alpha|(\xi)=\frac{r^n(\xi)\mathcal{K}(\vec{r}(\xi))}{h(\alpha(\xi))},
\end{align}
and
\begin{align}\label{2.2.8}
|Jac \alpha^*|(x)=\frac{h(x)}{r^n(\alpha^*(x))\mathcal{K}(\upsilon_M^{-1}(x))},
\end{align}

\section{\bf Gauss curvature flow and its associated functional}
In this section, we shall introduce an anisotropic Gauss-like curvature flow and its associated functional for solving the Orlicz-Aleksandrov problem.

First, we need to show that Orlicz-Aleksandrov problem is equivalent to solving a Monge-Amp\`{e}re type equation. From (\ref{2.1.3}), we know that the Orlicz-integral curvature $\mathcal{J}_\phi(K,\cdot)$ is absolutely continuous with respect to the integral curvature $\mathcal{J}(K,\cdot)$, namely,
\begin{align}\label{3.02}
d\mathcal{J}_\phi(K,\cdot)=\phi(\rho_K)d\mathcal{J}(K,\cdot).
\end{align}

If the convex body $K$ is $C^2$ smooth with positive Gauss curvature, then it follows from (\ref{2.1.2}) and (\ref{3.02}) that the $\mathcal{J}_\phi(K,\cdot)$ has a continuous density, given by
\begin{align}\label{3.03}
\frac{h_K\phi(1/h_K)}{(|\nabla h_K|^2+h_K^2)^{\frac{n}{2}}}\det(\nabla^2h_K+h_KI).
\end{align}
From (\ref{3.02}) and (\ref{3.03}), if the given measure $\mu$ on $S^{n-1}$ has a density $g$, then the equation (\ref{1.3}) is reduced into
\begin{align}\label{3.01}
\frac{\lambda h\phi(1/h)}{(|\nabla h|^2+h^2)^{\frac{n}{2}}}\det(\nabla^2h+hI)=g\ \ on \ \ S^{n-1}.
\end{align}

Let $M_0$ be a closed, smooth and strictly convex hypersurface in $\mathbb{R}^n$ enclosing the origin, and $\phi: (0,+\infty)\rightarrow(0,+\infty)$ be a positive smooth function. Consider the following anisotropic Gauss-like curvature flow
\begin{align}\label{3.0}
\begin{cases}
\frac{\partial X(x,t)}{\partial t}=-g(\nu)\frac{r^n}{\phi(r)}\mathcal{K}\nu+X(x,t),\\
X(x,0)=X_0(x).
\end{cases}
\end{align}

By the definition of support function, it is easy for us to know $h(x,t)=\langle x,X(x,t)\rangle$. From the evolution equation of $X(x,t)$ in (\ref{3.0}), we derive the  evolution equation of the corresponding support function $h(x,t)$
\begin{align}\label{2.3.1}
\frac{\partial h(x,t)}{\partial t}=-g(x)\frac{r^n}{\phi(r)}\mathcal{K}+h(x,t)\ \ on \ \ S^{n-1}\times[0,T).
\end{align}

Denote the radial function of $M_t$ by $\rho(u,t)$. For each $t$, let $u$ and $x$ be related through the following equality:
\begin{align}\label{2.3.2}
\rho(u,t)u=\nabla h(x,t)+h(x,t)x.
\end{align}
Thus, $x$ can be expressed as $x=x(u,t)$, by (\ref{2.3.2}), we get
\begin{align*}
\log \rho(u,t)=\log h(x,t)-\log\langle x,u\rangle.
\end{align*}
Differentiating the above identity, it is easy to see
\begin{align}\label{2.3.6}
\frac{1}{\rho(u,t)}\frac{\partial \rho(u,t)}{\partial t}=\frac{1}{h(x,t)}\frac{\partial h(x,t)}{\partial t}
\end{align}
Therefore, by the definitions of radial function and $r$, the normalised flow (\ref{3.0}) can be also described by the following scalar equation for $r(\cdot,t)$,
\begin{align}\label{2.3.7}
\frac{\partial r}{\partial t}(\xi,t)=-g(\xi)\frac{r^{n+1}}{h\phi(r)}\mathcal{K}+r(\xi,t)\ \ on \ \ S^{n-1}\times[0,T),
\end{align}
where $\mathcal{K}$ denotes the Gauss curvature at $r(\xi,t)\xi\in M_t$.

It can be checked that flow (\ref{3.0}) is the gradient flow of the functional given by
\begin{align*}
\mathcal{F}(M_t)=\int_{S^{n-1}}\log h(x,t)dx-\int_{S^{n-1}}\frac{\varphi(r(\xi,t))}{g(x)}d\xi,\ \ \ \ (\cdot,t)\in S^{n-1}\times[0,T),
\end{align*}
where
$$\varphi(t)=\int_0^t\frac{\phi(s)}{s}ds,$$
$h$ and $r$ are the support function and radial function of $M_t$ respectively.

\noindent{\bf Lemma 3.1}~~{\it Let $M_t$ be a strictly convex solution to the flow (\ref{3.0}).  Then the functional $\mathcal{F}(M_t)$ is non-increasing along the flow (\ref{3.0}). That is
\begin{align*}
\frac{\partial}{\partial t}\mathcal{F}(M_t)\leq0.
\end{align*}
with equalities if and only if the support function of $M_t$ satisfies the elliptic equation (\ref{3.01}).}
{\it \bf Proof.}~~From (\ref{2.3.1}) and (\ref{2.3.7}). By the fact that $r^nd\xi=\frac{h}{\mathcal{K}}dx$, we have
\begin{align*}
\frac{\partial}{\partial t}\mathcal{F}(M_t)&=\int_{S^{n-1}}\frac{1}{h(x,t)}\frac{\partial h}{\partial t}dx
-\int_{S^{n-1}}\frac{\phi(r(\xi,t))}{r(\xi,t)}\frac{\partial r}{\partial t}\frac{1}{g(x)}d\xi\\
&=\int_{S^{n-1}}\frac{1}{h}\frac{\partial h}{\partial t}dx
-\int_{S^{n-1}}\frac{\phi(r)}{g(x)}\frac{1}{r^n\mathcal{K}}\frac{\partial h}{\partial t}dx\\
&=\int_{S^{n-1}}\bigg(\frac{g(x)r^n\mathcal{K}-\phi(r)h}{g(x)hr^n\mathcal{K}}\bigg)\frac{\partial h}{\partial t}dx\\
&=\int_{S^{n-1}}\bigg(\frac{g(x)r^n\mathcal{K}/\phi(r)-h}{g(x)hr^n\mathcal{K}/\phi(r)}\bigg)\frac{\partial h}{\partial t}dx\\
&=-\int_{S^{n-1}}\frac{\bigg(\frac{g(x)r^n\mathcal{K}}{\phi(r)}-h\bigg)^2}{\frac{g(x)hr^n\mathcal{K}}{\phi(r)}}dx\\
&\leq0.
\end{align*}
Clearly $\frac{\partial}{\partial t}\mathcal{F}(M_t)=0$ holds if and only if
$$\frac{g(x)r^n\mathcal{K}}{\phi(r)}=h.$$
From (\ref{2.2.2}) and the concept of radial function, the above equation implies equation (\ref{3.01}) with $\lambda=1$. \hfill${\square}$

\section{\bf The long-time existence of the flow}
In this section, we will obtain the long-time existence of the flow (\ref{3.0}). It is equivalent to obtain the long-time existence of the evolution equation (\ref{2.3.1}). The main work is to obtain the $C^0$, $C^1$ and $C^2$-estimates for the (\ref{2.3.1}).

\noindent{\bf 4.1~~$C^0, C^1$-Estimates}

The following lemma obtains the $C^0$-estimate.

\noindent{\bf Lemma 4.1}~~{\it Let $h$ be a smooth solution of (\ref{2.3.1}), and $g$ be a positive, smooth function on $S^{n-1}$ satisfying (\ref{1.5}), then there is a
positive constant $C$ independent of $t$ such that
\begin{align}\label{3.1.1}
\frac{1}{C}\leq h(x,t)\leq C,
\end{align}
and
\begin{align}\label{3.1.2}
\frac{1}{C}\leq r(\xi,t)\leq C
\end{align}
for $\forall(\cdot,t)\in S^{n-1}\times(0,T]$.}

\noindent{\it \bf Proof.}~~Since (\ref{3.1.1}) and (\ref{3.1.2}) are equivalent, hence,
for upper bound (or lower bound) we only need to establish (\ref{3.1.1}) or (\ref{3.1.2}). Suppose that $h(x,t)$ is maximized at point $x_1\in S^{n-1}$, then at $x_1$, we get
\begin{align*}
\nabla h=0,\ \ \ \ \nabla^2h\leq0\ \ \ \ and \ \ \ \ r=h.
\end{align*}
From (\ref{2.3.1}), at $x_1$, we have
\begin{align*}
\frac{\partial h}{\partial t}=-g(x)\mathcal{K}\frac{h^{n}}{\phi(h)}+h
 \leq-g(x)\frac{h}{\phi(h)}+h
=\frac{h}{\phi(h)}\bigg(\phi(h)-g(x)\bigg).
\end{align*}
Taking $\overline{\Lambda}=\limsup_{s\rightarrow+\infty}\varphi(s)$. By (\ref{1.5}), $\varepsilon=\frac{1}{2}(\min_{S^{n-1}}g(x)-\overline{\Lambda})$ is positive and there exists some positive constant $C_1>0$ such that
\begin{align*}
\phi(h)<\overline{\Lambda}+\varepsilon
\end{align*}
for $h<C_1$. This together with (\ref{1.5})
\begin{align*}
\phi(h)-g(x)<\overline{\Lambda}+\varepsilon-\min_{S^{n-1}}g(x)<0,
\end{align*}
which implies that at maximal point
\begin{align*}
\frac{\partial h}{\partial t}<0
\end{align*}
Therefore
\begin{align*}
h\leq\max \{C_1,\max_{S^{n-1}}h(x,0)\}.
\end{align*}

Similarly, we can estimate $\min_{S^{n-1}}h(x,t)$. Suppose that $h(x,t)$ is minimized at point $x_2\in S^{n-1}$, then at $x_2$, we get
\begin{align*}
\nabla h=0,\ \ \ \ \nabla^2h\geq0\ \ \ \ and \ \ \ \ r=h.
\end{align*}
From (\ref{2.3.1}), at $x_2$, we have
\begin{align*}
\frac{\partial h}{\partial t}=-g(x)\mathcal{K}\frac{h^{n}}{\phi(h)}+h
 \geq-g(x)\frac{h}{\phi(h)}+h
=\frac{h}{\phi(h)}\bigg(\phi(h)-g(x)\bigg).
\end{align*}
Taking $\underline{\Lambda}=\liminf_{s\rightarrow0^+}\varphi(s)$. By (\ref{1.5}), $\varepsilon=\frac{1}{2}(\underline{\Lambda}-\max_{S^{n-1}}g(x))$ is positive and there exists some positive constant $C_2>0$ such that
\begin{align*}
\phi(h)>\underline{\Lambda}-\varepsilon
\end{align*}
for $h<C_2$. This together with (\ref{1.5})
\begin{align*}
\phi(h)-g(x)>\underline{\Lambda}-\varepsilon-\max_{S^{n-1}}g(x)>0,
\end{align*}
which implies that at minimal point
\begin{align*}
\frac{\partial h}{\partial t}>0
\end{align*}
Therefore
\begin{align*}
h\geq\min \{C_2,\min_{S^{n-1}}h(x,0)\}.
\end{align*}
The proof of the lemma is completed.      \hfill${\square}$

Since the convexity of $M_t$, combining with Lemma4.1, we can obtain the $C^1$-estimates as follows.

\noindent{\bf Lemma 4.2}~~{\it Under the assumption of Lemma 4.1, we have
\begin{align*}
|\nabla h(x,t)|\leq C,\ \ \ \ and \ \ \ \ |\nabla r(\xi,t)|\leq C,
\end{align*}
for $\forall(\cdot,t)\in S^{n-1}\times (0,T]$. Here $C$ is a positive constant depending only on the constant in Lemma 4.1.}

\noindent{\it \bf Proof.}~~By virtue of the (\ref{2.2.1}), (\ref{2.2.11}) and (\ref{2.2.2}), there is
$$r^2=h^2+|\nabla h|^2,$$
which implies that
$$|\nabla h|\leq r.$$

Then from (\ref{2.2.3}), we have
$$h=\frac{r^2}{\sqrt{r^2+|\nabla r|^2}},$$
which implies that
$$|\nabla r|=\frac{|\nabla h|}{h}r\leq\frac{r^2}{h}.$$
From the Lemma 4.1, we directly obtain the estimate of this lemma.  \hfill${\square}$

\noindent{\bf 4.2~~$C^2$-Estimate}

In this subsection, we will establish the upper and lower bound of principal curvatures. This estimates can be obtained by considering proper auxiliary functions; see, e.g., \cite{CH,LL} for similar techniques. We take a local orthonormal frame $\{e_1,...,e_{n-1}\}$ on $S^{n-1}$ such that the standard metric on $S^{n-1}$ is $\{\delta_{ij}\}$. We first derive an upper bound for the Gauss curvature.

\noindent{\bf Lemma 4.3}~~{\it Let $h$ be a smooth solution of (\ref{2.3.1}), and $g$ be a positive, smooth function on $S^{n-1}$ satisfying (\ref{1.5}), then there is a
positive constant $C$ independent of $t$ such that
\begin{align*}
\mathcal{K}(x,t)\leq C,
\end{align*}
for $\forall (x,t)\in S^{n-1}\times[0,T)$.}

\noindent{\it \bf Proof.}~~Let us consider the following auxiliary function
\begin{align}\label{3.1}
\Theta(x,t)=\frac{-\partial_th+h}{h-\varepsilon_0}=g(x)\frac{r^n}{\phi(r)}\frac{\mathcal{K}}{h-\varepsilon_0}
\end{align}
where $\varepsilon_0$ is a positive constant satisfying $\varepsilon_0=\frac{1}{2}\inf h(x,t)$, $\forall(x,t)\in S^{n-1}\times[0,T)$.

From (\ref{3.1}), the upper bound of $\mathcal{K}$ follows from $\Theta(x,t)$. Hence we only need to derive the upper bound of $\Theta(x,t)$. At any maximum of $\Theta$ at $x_0$ we have
\begin{align}\label{3.2}
0=\nabla_i\Theta=\frac{-\partial_th_i+h_i}{h-\varepsilon_0}+\frac{(\partial_th-h)h_i}{(h-\varepsilon_0)^2},
\end{align}
and using the (\ref{3.2}), then
\begin{align}\label{3.3}
0\geq\nabla_{ij}\Theta=\frac{-\partial_th_{ij}+h_{ij}}{h-\varepsilon_0}+\frac{(\partial_th-h)h_{ij}}{(h-\varepsilon_0)^2},
\end{align}
where $\nabla_{ij}\Theta\leq0$ should be understood in the sense of negative semi-definite matrix. As in the background metrial, we know the fact $b_{ij}=h_{ij}+h\delta_{ij}$, and $b^{ij}$ its inverse matrix, which together with the (\ref{3.3}), we can get
\begin{align*}
\partial_tb_{ij}&=\partial_th_{ij}+\partial_th\delta_{ij}\\
&\geq h_{ij}+\frac{(\partial_th-h)h_{ij}}{h-\varepsilon_0}+\partial_th\delta_{ij}\\
&=b_{ij}+\frac{(\partial_th-h)h_{ij}}{h-\varepsilon_0}+(\partial_th-h)\delta_{ij}\\
&=b_{ij}+\frac{\partial_th-h}{h-\varepsilon_0}(h_{ij}+h\delta_{ij}-\varepsilon_0\delta_{ij})\\
&=b_{ij}-\Theta(b_{ij}-\varepsilon_0\delta_{ij}).
\end{align*}
By the fact (\ref{2.2.6}), we obtain
\begin{align}\label{3.4}
\nonumber \partial_t\mathcal{K}&=-\mathcal{K}b^{ij}\partial_tb_{ij}\\
&\leq-\mathcal{K}b^{ij}[b_{ij}-\Theta(b_{ij}-\varepsilon_0\delta_{ij})]\\
\nonumber &=-\mathcal{K}[(n-1)(1-\Theta)+\Theta\varepsilon_0\mathcal{H}],
\end{align}
where $\mathcal{H}$ denotes the mean curvature of $X(\cdot,t)$.

From the (\ref{3.1}) and Lemma 4.1,  there exists a constant $C_1$ such that
\begin{align}\label{3.5}
\frac{1}{C_1}\Theta(x,t)\leq \mathcal{K}(x,t)\leq C_1\Theta(x,t),
\end{align}
where $C_1$ is a positive constant. Noting
\begin{align*}
\frac{1}{n-1}\mathcal{H}\geq \mathcal{K}^{\frac{1}{n-1}},
\end{align*}
and combining the inequalities (\ref{3.4}), we obtain
\begin{align}\label{3.6}
\partial_t\mathcal{K}&\leq(n-1)\mathcal{K}\Theta-(n-1)\varepsilon_0\Theta \mathcal{K}^{\frac{n}{n-1}}.
\end{align}

Now we estimate $\partial_t\Theta$. From (\ref{3.1}), we have
\begin{align}\label{3.7}
\partial_t\Theta=\partial_t\bigg(\frac{r^n}{(h-\varepsilon_0)\phi(r)}\bigg)g(x)\mathcal{K}
+\frac{g(x)r^n}{(h-\varepsilon_0)\phi(r)}\partial_t\mathcal{K},
\end{align}
where
$$\partial_th=h-(h-\varepsilon_0)\Theta,$$
\begin{align*}
\partial_t\bigg(\frac{r^n}{(h-\varepsilon_0)\phi(r)}\bigg)
=\frac{nr^{n-1}}{(h-\varepsilon_0)\phi(r)}\partial_tr-\frac{(h-\varepsilon_0)r^n\phi^\prime(r)}{((h-\varepsilon_0)\phi(r))^2}\partial_tr
-\frac{r^n\phi(r)}{((h-\varepsilon_0)\phi(r))^2}\partial_th.
\end{align*}
From (\ref{3.6}) and (\ref{3.5}), we have at $x_0$
\begin{align*}
\partial_t\Theta&\leq(h-\varepsilon_0)\bigg(\frac{r\phi^\prime(r)}{h\phi(r)}-\frac{n}{h}+\frac{1}{h-\varepsilon_0}\bigg)\Theta^2
+\frac{g(x)r^n}{(h-\varepsilon_0)\phi(r)}\bigg((n-1)\mathcal{K}\Theta-(n-1)\varepsilon_0\Theta \mathcal{K}^{\frac{n}{n-1}}\bigg)\\
&\leq C_2\Theta^2+\frac{g(x)r^n}{(h-\varepsilon_0)\phi(r)}\bigg((n-1)\mathcal{K}\Theta-(n-1)\varepsilon_0\Theta \mathcal{K}^{\frac{n}{n-1}}\bigg)\\
&\leq C_3\Theta^2\bigg(C_4-\varepsilon_0\Theta^{\frac{1}{n-1}}\bigg),
\end{align*}
where $C_2, C_3$ and $C_4$ are positive constant depending only on the constant $C$ in Lemma 4.1, and the upper and lower bounds of $g$ on $S^{n-1}$ and $\phi$ on $[1/C,C]$.

Now one can see that whenever $\Theta>\bigg(\frac{C_4}{\varepsilon_0}\bigg)^{n-1}$ which is independent of $t$,
\begin{align*}
\partial_t\Theta<0,
\end{align*}
which implies that $\Theta$ has  a uniform upper bound.

For any $(x,t)$
$$\mathcal{K}(x,t)=\frac{(h-\varepsilon_0)\phi(r)\Theta(x,t)}{g(x)r^n}\leq\frac{(h-\varepsilon_0)\phi(r)\Theta(x_0,t)}{g(x)r^n}\leq C,$$
namely, $\mathcal{K}$ has a uniform upper bound.     \hfill${\square}$

Now, we estimate the principal curvatures are bounded from below along the flow (\ref{3.0}). To obtain the positive lower bound for the principal curvatures of $M_t$, we will study an expanding flow by Gauss curvature for the dual hypersurface of $M_t$.

\noindent{\bf Lemma 4.4}~~{\it Under the conditions of Lemma 4.3, then the principal curvature $k_i$ for $i=1,...,n-1$ satisfies
\begin{align*}
\frac{1}{C}\leq k_i\leq C,
\end{align*}
where $C$ is a positive constant independent of $t$.}

\noindent{\it \bf Proof.}~~To prove the lower bound of $k_i$, we employ the dual flow of (\ref{3.0}), and establish an upper bound of principal curvature for the dual flow. This together with Lemma 4.3, also implies the upper bound of $k_i$.

We denote by $M_t^*$ the polar set of $M_t=X(S^{n-1},t)$. From the definition of polar set, if $r(\cdot,t)$ is the radial function of $M_t$, then
\begin{align}\label{3.9}
r(\xi,t)=\frac{1}{h^*(\xi,t)},
\end{align}
where $h^*(\xi,t)$ denotes the supper function of $M_t^*$. it is well-know that $\alpha_{M_t^*}=\alpha_{M_t}^*$, see e.g. \cite{H}. This implies $|Jac \alpha_{M_t}^*||Jac \alpha_{M_t}|=1$, thus by (\ref{2.2.7}) and (\ref{2.2.8}) we have, under a local orthonormal frame on $S^{n-1}$
\begin{align}\label{3.10}
\frac{(h^*(\xi,t))^{n+1}h^{n+1}(x,t)}{\mathcal{K}^*(p^*)\mathcal{K}(p)}=1,
\end{align}
where $p\in M_t$, $p^*\in M_t^*$ are two points satisfying $p\cdot p^*=1$, and $x, \xi$ are the unit outer normals of $M_t$ and $M_t^*$ at $p$ and $p^*$. Therefore by equation (\ref{2.3.7}), we obtain the equation for $h^*$,
\begin{align}\label{3.11}
\partial_th^*(\xi,t)=g(\xi)\frac{(h^*(\xi,t))^2}{\phi(r^*)(r^*)^{n}\mathcal{K}^*}-h^*(\xi,t),
\ \ \ \ \forall(\cdot,t)\in S^{n-1}\times(0,T],
\end{align}
where $\mathcal{K}^*=(\det(\nabla^2h^*+h^*I))^{-1}$ is the Gauss curvature of $M_t^*$ at the point $p^*=\nabla h^*(\xi,t)+h^*(\xi,t)\xi$, and
\begin{align*}
r^*=|p^*|=\sqrt{|\nabla h^*|^2+(h^*)^2}(\xi,t)
\end{align*}
is the distance from $p^*$ to the origin. Note that $g$ takes value at
\begin{align*}
x=\frac{p^*}{|p^*|}=\frac{\nabla h^*+h^*\xi}{\sqrt{|\nabla h^*|^2+(h^*)^2}}\in S^{n-1}.
\end{align*}

By (\ref{3.9}), $\frac{1}{C}\leq h^*\leq C$ and $|\nabla h^*|\leq C$ for some $C$ only depending on $\max_{S^{n-1}\times(0,T]}h$ and $\min_{S^{n-1}\times(0,T]}h$.

Let $b_{ij}^*=h_{ij}^*+h^*\delta_{ij}$, and $b_*^{ij}$ be the inverse matrix of $b_{ij}^*$. As discussed in Section 2, the eigenvalues of $b_{ij}^*$ and $b_*^{ij}$ are respectively the principal radii and principal curvature of $M_t^*$. Consider the following function
\begin{align}\label{3.12}
W(\xi,t,\tau)=\log b_*^{\tau\tau}-\beta\log h^*+\frac{A}{2}(r^*)^2,
\end{align}
where $\tau$ is a unit vector in the tangential space of $S^{n-1}$, while $\beta$ and $A=A(\beta)$ are large constants to be specified later on. Assume $w$ attain its maximum at $(\xi_0,t_0)$, along the direction $\tau=e_1$. By rotation, we also assume $b_*^{ij}$ and $b_{ij}^*$ are diagonal at this point.

It is direct to see, at the point where $W$ attains its maximum
\begin{align}\label{3.13}
0=\nabla_iW=-b_*^{11}b_{11;i}^*-\beta\frac{h_i^*}{h^*}+Ar^*r_i^*,
\end{align}
and
\begin{align}\label{3.14}
0\geq\nabla_{ij}W=-b_*^{11}b_{11;ij}^*-(b_*^{11})^2b_{11;i}^*b_{11;j}^*-\beta\bigg(\frac{h_{ij}^*}{h^*}-\frac{h_i^*h_j^*}{(h^*)^2}\bigg)+
A(r^*r_{ij}^*+r_i^*r_j^*),
\end{align}
where $\nabla_{ij}W\leq0$ should be understood in the sence of negative semi-definite matrix. Note that $b_{ij;k}^*$ is symmetric in all indices.  Without loss of generality, if we assume $t_0>0$, then at $(\xi_0,t_0)$, we also have
\begin{align}\label{3.15}
\nonumber0\leq\partial_tW&=b_{11}^*\partial_tb_*^{11}-\beta\frac{\partial_th^*}{h^*}+Ar^*\partial_tr^*\\
&=-b_*^{11}\partial_tb_{11}^*-\beta\frac{\partial_th^*}{h^*}+Ar^*\partial_tr^*.
\end{align}

We can rewrite the flow (\ref{3.11}) as
\begin{align}\label{3.16}
\log(\partial_th^*+h^*)=\log S_n+\alpha(\xi,t),
\end{align}
where
\begin{align*}
\alpha(\xi,t)=\log\bigg(g(\xi)\frac{(h^*)^2}{\phi(r^*)(r^*)^{n}}\bigg).
\end{align*}
Differentiating (\ref{3.16}) gives
\begin{align}\label{3.17}
\frac{\partial_th_k^*+h_k^*}{\partial_th^*+h^*}=b_*^{ij}b_{ij;k}^*+\nabla_k\alpha,
\end{align}
and
\begin{align}\label{3.18}
\frac{\partial_th_{11}^*+h_{11}^*}{\partial_th^*+h^*}=\frac{(\partial_th_1^*+h_1^*)^2}{(\partial_th^*+h^*)^2}
+b_*^{ij}b_{ij;11}^*-b_*^{ii}b_*^{jj}(b_{ij,1}^*)^2+\nabla_{11}\alpha.
\end{align}

Dividing (\ref{3.15}) by $\partial_th^*+h^*$ and using (\ref{3.18}), we have
\begin{align}\label{3.19}
\nonumber0&\leq-b_*^{11}\bigg(\frac{\partial_th_{11}^*+h_{11}^*}{\partial_th^*+h^*}-\frac{b_{11}^*}{\partial_th^*+h^*}+1\bigg)-
\frac{\beta\partial_th^*}{h^*(\partial_th^*+h^*)}+A\frac{r^*\partial_tr^*}{\partial_th^*+h^*}\\
&=-b_*^{11}\frac{\partial_th_{11}^*+h_{11}^*}{\partial_th^*+h^*}-b_*^{11}+\frac{1+\beta}{\partial_th^*+h^*}
-\frac{\beta}{h^*}+A\frac{r^*\partial_tr^*}{\partial_th^*+h^*}\\
\nonumber&\leq-b_*^{11}b_*^{ij}b_{ij;11}^*+b_*^{11}b_*^{ii}b_*^{jj}(b_{ij;1}^*)^2-b_*^{11}\nabla_{11}\alpha
+\frac{1+\beta}{\partial_th^*+h^*}+A\frac{r^*\partial_tr^*}{\partial_th^*+h^*}.
\end{align}
By the Ricci identity, we have
\begin{align*}
b_{ij;11}^*=b_{11;ij}^*-\delta_{ij}b_{11}^*+\delta_{11}b_{ij}^*-\delta_{i1}b_{1j}^*+\delta_{1j}b_{1i}^*.
\end{align*}
Plugging this identity in (\ref{3.19}), and employing (\ref{3.14}), we obtain
\begin{align}\label{3.20}
\nonumber0&\leq b_{*}^{11}\bigg(b_*^{11}b_*^{ii}(b_{11;i}^*)^2-b_*^{ii}b_*^{jj}(b_{ij;1}^*)^2\bigg)+(\mathcal{H}^*-(n-1)b_*^{11})\\
\nonumber&\ \ \ \ -\beta \mathcal{H}^*+C\beta-\beta b_*^{ij}\frac{h_i^*h_j^*}{(h^*)^2}-b_*^{11}\nabla_{11}\alpha+\frac{1+\beta}{\partial_th^*+h^*}\\
&\ \ \ \ +A\frac{r^*\partial_tr^*}{\partial_th^*+h^*}-Ab_*^{ij}(r^*r^*_{ij}+r^*_ir^*_j)\\
\nonumber&\leq-\beta \mathcal{H}^*+C\beta-b_*^{11}\nabla_{11}\alpha+\frac{1+\beta}{\partial_th^*+h^*}
+A\frac{r^*\partial_tr^*}{\partial_th^*+h^*}-Ab_*^{ij}(r^*r^*_{ij}+r^*_ir^*_j),
\end{align}
where $\mathcal{H}^*=\sum b_*^{ii}$ is the mean curvature of $M_t^*$.

It is direct to calculate
$$r^*_t=\frac{h^*\partial_th^*+\sum h_k^*\partial_th_{k}^*}{r^*},$$
\begin{align}\label{3.21}
r^*_i=\frac{h^*h_i^*+\sum h_k^*h_{ki}^*}{r^*}=\frac{h_i^*b_{ii}^*}{r^*},
\end{align}
$$r^*_{ij}=\frac{h^*h_{ij}^*+h_i^*h_j^*+\sum h_k^*h_{kij}+\sum h_{ki}^*h_{kj}^*}{r^*}-\frac{h_i^*h_j^*b_{ii}^*b_{jj}^*}{(r^*)^3}.$$
Hence, by (\ref{3.17})
\begin{align*}
\frac{r^*r_t^*}{h_t^*+h^*}-b_*^{ij}(r^*r^*_{ij}+r_i^*r_j^*)&=\frac{h^*\partial_th^*}{\partial_th^*+h^*}-h^*b_*^{ij}h_{ij}^*-b_*^{ii}(h_{ii}^*)^2\\
&\ \ \ \ -\frac{|\nabla h^*|^2}{\partial_th^*+h^*}+\sum h_k^*\nabla_k\alpha.
\end{align*}
Since
\begin{align*}
\frac{h^*\partial_th^*}{\partial_th^*+h^*}-\frac{|\nabla h^*|^2}{\partial_th^*+h^*}=h^*-\frac{(r^*)^2}{\partial_th^*+h^*},
\end{align*}
and
\begin{align*}
-h^*b_*^{ij}h_{ij}^*-b_*^{ii}(h_{ii}^*)^2&=-h^*b_*^{ii}(b_{ii}^*-h^*\delta_{ii})-b_*^{ii}(b_{ii}^*-h^*\delta_{ii})^2\\
&=(n-1)h^*-\sum b_{ii}^*,
\end{align*}
we further deduce
\begin{align}\label{3.22}
\frac{r^*\partial_tr^*}{\partial_th^*+h^*}-b_*^{ij}(r^*r^*_{ij}+r^*_ir^*_j)\leq C-\frac{(r^*)^2}{\partial_th^*+h^*}+\sum h_k^*\nabla_k\alpha.
\end{align}
Plugging (\ref{3.22}) in (\ref{3.20}), we get
\begin{align}\label{3.23}
\nonumber0&\leq-\beta \mathcal{H}^*+C\beta+CA-b_*^{11}\nabla_{11}\alpha+\frac{1+\beta-A(r^*)^2}{\partial_th^*+h^*}+A\sum h_k^*\nabla_k\alpha\\
&\leq-\beta \mathcal{H}^*+C\beta+CA-b_*^{11}\nabla_{11}\alpha+A\sum h_k^*\nabla_k\alpha,
\end{align}
provided $A>2(1+\beta)/\min_{S^{n-1}\times(0,T]}(r^*)^2\geq C(1+\beta)$ for some $C>0$ only depending on $\max_{S^{n-1}\times(0,T]}h$.

By (\ref{3.13}) and (\ref{3.21}), we have
\begin{align*}
-b_*^{11}\nabla_{11}\alpha+A\sum h_k^*\nabla_k\alpha&\leq Cb_*^{11}(1+(h_{11}^*)^2)+CA-b_*^{11}\sum\alpha_{h_k^*}h_{k11}^*
+A\sum\alpha_{h_k^*}h_k^*h_{kk}^*\\
&\leq Cb_*^{11}+\frac{C}{b_*^{11}}+CA+C\beta.
\end{align*}
Hence (\ref{3.23}) can be further estimated as
\begin{align*}
0&\leq-\beta \mathcal{H}^*+Cb_*^{11}+C\beta+CA\\
&\leq(-\beta+C)b_*^{11}+C\beta+CA,
\end{align*}
by choosing $\beta$ large. This inequality tells us the principal curvature of $M^*$ are bounded from above, namely
\begin{align*}
\max_{\xi\in S^{n-1}}k_i^*(\xi,t)\leq C,\ \ \ \ \forall t\in(0,T],\ \ i=1,2,...,n-1.
\end{align*}
By the Lemma 4.3 and (\ref{3.10}), we have $\mathcal{K}^*(\cdot,t)\geq\frac{1}{C}$. Therefore
\begin{align*}
\frac{1}{C}\leq k_i^*(\cdot,t)\leq C,\ \ \ \ \forall (\cdot,t)\in S^{n-1}\times(0,T],\ \ i=1,2,...,n-1.
\end{align*}
By duality, Lemma 4.4 follows.  \hfill${\square}$

As a consequence of the above a priori estimates, one sees that the convexity of the hypersurface $M_t$ is preserved under the flow (\ref{3.0}) and the solution is uniformly convex.

Now we have proved that the principal curvatures of $M_t$ have uniform positive upper and lower bounds, this together with Lemmas 4.1 and 4.2 implies that the evolution equation (\ref{2.3.1}) is uniformly parabolic on any finite time interval. Thus, the result of \cite{K} and the standard parabolic theory show that the smooth solution of (\ref{2.3.1}) exists for all time, namely, flow (\ref{3.0}) has a long-time solution. And by these estimates again, a subsequence of $M_t$ converges in $C^\infty$ to a positive, smooth, uniformly convex hypersurface $M^\infty$ in $\mathbb{R}^n$.

\section{\bf Existence of the solutions to the Monge-Amp\`{e}re equation}
~~~~In this section, we complete proof of Theorem 1.2, namely we will prove the support function $\widetilde{h}$ of $M_\infty$ satisfies the following Monge-Amp\`{e}re equation:
\begin{align}\label{4.1}
\frac{h\phi(1/h)}{(|\nabla h|^2+h^2)^{\frac{n}{2}}}\det(\nabla^2h+hI)=g,
\end{align}

Recalling the functional $\mathcal{F}(M_t)$ we defined in Section 3
\begin{align*}
\mathcal{F}(M_t)=\int_{S^{n-1}}\log h(x,t)dx-\int_{S^{n-1}}\frac{\varphi(r(\xi,t))}{g(x)}d\xi,\ \ \ \ t\in[0,T).
\end{align*}
From the Lemma 3.1, there exists a positive constant $C$ which is independent of $t$, such that
\begin{align}\label{4.2}
\mathcal{F}(M_t)\leq C.
\end{align}

Since $\mathcal{F}(M_t)$ is non-increasing for any $t>0$. From
\begin{align*}
\int_0^t\bigg(-\frac{\partial}{\partial t}\mathcal{F}(M_t)\bigg)dt=\mathcal{F}(M_0)-\mathcal{F}(M_t)\leq\mathcal{F}(M_0),
\end{align*}
we have
\begin{align*}
\int_0^\infty\bigg(-\frac{\partial}{\partial t}\mathcal{F}(M_t)\bigg)dt\leq\mathcal{F}(M_0),
\end{align*}
this implies that there exists a subsequence of times $t_j\rightarrow\infty$ such that
\begin{align*}
-\frac{\partial}{\partial t}\mathcal{F}(M_{t_j})\rightarrow0\ \ \ \ as\ \ \ \ t_j\rightarrow\infty.
\end{align*}
Recalling Lemma 3.1
$$\frac{\partial \mathcal{F}(M_t)}{\partial t}
=-\int_{S^{n-1}}\frac{\bigg(\frac{g(x)r^n\mathcal{K}}{\phi(r)}-h\bigg)^2}{\frac{g(x)hr^n\mathcal{K}}{\phi(r)}}dx.$$
Since $h, r$ and $\mathcal{K}$ have uniform positive upper and lower bounds, by passing to the limit, we obtain
$$\frac{g(x)r^n\widetilde{\mathcal{K}}}{\phi(r)}=\widetilde{h},$$
where $\widetilde{h}$ and $\widetilde{\mathcal{K}}$ are the support function and Gauss curvature of $M_\infty$. Namely
$$g(x)\frac{\bigg(\sqrt{|\nabla \widetilde{h}|^2+\widetilde{h}^2}\bigg)^n}{\phi(1/\widetilde{h})}\widetilde{\mathcal{K}}=\widetilde{h}\ \ \ \ on\ \ \ \ S^{n-1,}$$
which is just equation (\ref{4.1}). The proof of Theorem 1.2 is now completed.

At the same time, for Theorem 1.1, we have showed that for smooth $\phi$ and $g$, there exists a smooth solution $h$ to (\ref{1.3}) with $\lambda=1$.

Finally, we provide a special uniqueness of the Orlicz-Aleksandrov problem under the appropriate condition.

\noindent{\bf Theorem 5.1}~~{\it Assume $\phi$ is a positive, continuous function. If whenever
\begin{align}\label{4.51}
\phi(cs^{-1})\leq \phi(s^{-1})
\end{align}
hold for some positive $c, s$, there must be $c\geq1$. Then the solution to the
\begin{align}\label{4.5}
\frac{h\phi(1/h)}{(|\nabla h|^2+h^2)^{\frac{n}{2}}}\det(\nabla^2h+hI)=g.
\end{align}
is unique.}

\noindent{\it \bf Proof of Theorem 5.1.}~~Let $h_1$ and $h_2$ are two solutions of (\ref{4.5}).
Assume $\frac{h_1}{h_2}$ attain its maximum at $x_0\in S^{n-1}$. Taking $Q=\log\frac{h_1}{h_2}$, then at $x_0$
\begin{align*}
0=\nabla Q=\frac{\nabla h_1}{h_1}-\frac{\nabla h_2}{h_2},
\end{align*}
and
\begin{align*}
0\geq\nabla^2 Q=\frac{\nabla^2 h_1}{h_1}-\frac{\nabla^2 h_2}{h_2}.
\end{align*}
By the equation (\ref{4.5}), we have at $x_0$
\begin{align}\label{4.7}
\nonumber1&=\frac{\det(\nabla^2h_2+h_2I)\bigg(|\nabla h_2|^2+h_2^2\bigg)^{-\frac{n}{2}}\phi(h_2^{-1})h_2}{\det(\nabla^2h_1+h_1I)\bigg(|\nabla h_1|^2+h_1^2\bigg)^{-\frac{n}{2}}\phi(h_1^{-1})h_1}\\
&=\frac{h_2^{n}\det(\frac{\nabla^2 h_2}{h_2}+I)\bigg[h_2^2\bigg(|\frac{\nabla h_2}{h_2}|^2+1\bigg)\bigg]^{-\frac{n}{2}}\phi(h_2^{-1})}{h_1^{n}\det(\frac{\nabla^2 h_1}{h_1}+I)\bigg[h_1^2\bigg(|\frac{\nabla h_1}{h_1}|^2+1\bigg)\bigg]^{-\frac{n}{2}}\phi(h_1^{-1})}\\
\nonumber&\geq\frac{\phi(h_2^{-1})}{\phi(h_1^{-1})}.
\end{align}
Write $h_2(x_0)=ch_1(x_0)$, then the above inequality reads
$$\phi(h_1^{-1})\geq\phi(ch_1^{-1}).$$
By our assumption (\ref{4.51}), we have $c\geq1$. Namely $h_1(x_0)\geq h_2(x_0)$.

Interchanging $h_1$ and $h_2$, then $h_2(x_0)\geq h_1(x_0)$. Therefore, we have $h_1\equiv h_2$.         \hfill${\square}$

\end{document}